\documentclass[11pt]{amsart}
\usepackage{amssymb,mathrsfs}
\title[\tiny {Theta Sums of Higher Index }]
{Theta Sums of Higher Index}

\author{Jae-Hyun Yang}
\address{Department of Mathematics, Inha University, Incheon
22212, Korea}

\email{jhyang@inha.ac.kr }

\begin{document}

\newtheorem{theorem}{Theorem}[section]
\newtheorem{corollary}{Corollary}[section]
\newtheorem{lemma}{Lemma}[section]
\newtheorem{proposition}{Proposition}[section]
\newtheorem{remark}{Remark}[section]
\newtheorem{definition}{Definition}[section]

\renewcommand{\theequation}{\thesection.\arabic{equation}}
\renewcommand{\thetheorem}{\thesection.\arabic{theorem}}
\renewcommand{\thelemma}{\thesection.\arabic{lemma}}
\newcommand{\BR}{\mathbb R}
\newcommand{\BQ}{\mathbb Q}
\newcommand{\bn}{\bf n}
\def\charf {\mbox{{\text 1}\kern-.24em {\text l}}}
\newcommand{\BC}{\mathbb C}
\newcommand{\BZ}{\mathbb Z}

\thanks{
\\ \indent
The author was supported by Basic Science Program through the National Research Foundation of Korea
\\ \indent
(NRF) funded by the Ministry of Education, Science
and Technology\,(49562-1) and also by INHA UNI-\\
 \indent VERSITY Research Grant.}


\begin{abstract}
{In this paper, we obtain some behaviours of theta sums of higher index for the Schr{\"o}dinger-Weil representation
of the Jacobi group associated with a positive definite symmetric real matrix of degree $m$.}
\end{abstract}
\maketitle

\newcommand\tr{\triangleright}
\newcommand\al{\alpha}
\newcommand\be{\beta}
\newcommand\g{\gamma}
\newcommand\gh{\Cal G^J}
\newcommand\G{\Gamma}
\newcommand\de{\delta}
\newcommand\e{\epsilon}
\newcommand\z{\zeta}
\newcommand\vth{\vartheta}
\newcommand\vp{\varphi}
\newcommand\om{\omega}
\newcommand\p{\pi}
\newcommand\la{\lambda}
\newcommand\lb{\lbrace}
\newcommand\lk{\lbrack}
\newcommand\rb{\rbrace}
\newcommand\rk{\rbrack}
\newcommand\s{\sigma}
\newcommand\w{\wedge}
\newcommand\fgj{{\frak g}^J}
\newcommand\lrt{\longrightarrow}
\newcommand\lmt{\longmapsto}
\newcommand\lmk{(\lambda,\mu,\kappa)}
\newcommand\Om{\Omega}
\newcommand\ka{\kappa}
\newcommand\ba{\backslash}
\newcommand\ph{\phi}
\newcommand\M{{\Cal M}}
\newcommand\bA{\bold A}
\newcommand\bH{\bold H}

\newcommand\Hom{\text{Hom}}
\newcommand\cP{\Cal P}
\newcommand\cH{\Cal H}

\newcommand\pa{\partial}

\newcommand\pis{\pi i \sigma}
\newcommand\sd{\,\,{\vartriangleright}\kern -1.0ex{<}\,}
\newcommand\wt{\widetilde}
\newcommand\fg{\frak g}
\newcommand\fk{\frak k}
\newcommand\fp{\frak p}
\newcommand\fs{\frak s}
\newcommand\fh{\frak h}
\newcommand\Cal{\mathcal}

\newcommand\fn{{\frak n}}
\newcommand\fa{{\frak a}}
\newcommand\fm{{\frak m}}
\newcommand\fq{{\frak q}}
\newcommand\CP{{\mathcal P}_g}
\newcommand\Hgh{{\mathbb H}_g \times {\mathbb C}^{(h,g)}}
\newcommand\BD{\mathbb D}
\newcommand\BH{\mathbb H}
\newcommand\CCF{{\mathcal F}_g}
\newcommand\CM{{\mathcal M}}
\newcommand\Ggh{\Gamma_{g,h}}
\newcommand\Chg{{\mathbb C}^{(h,g)}}
\newcommand\Yd{{{\partial}\over {\partial Y}}}
\newcommand\Vd{{{\partial}\over {\partial V}}}

\newcommand\Ys{Y^{\ast}}
\newcommand\Vs{V^{\ast}}
\newcommand\LO{L_{\Omega}}
\newcommand\fac{{\frak a}_{\mathbb C}^{\ast}}

\renewcommand\th{\theta}
\renewcommand\l{\lambda}
\renewcommand\k{\kappa}
\newcommand\tg{\widetilde\gamma}
\newcommand\wmo{{\mathscr W}_{\mathcal M,\Omega}}
\newcommand\hrnm{H_\BR^{(n,m)}}
\newcommand\rmn{\BR^{(m,n)}}

\begin{section}{{\bf Introduction}}
\setcounter{equation}{0}

For a given fixed positive integer $n$, we let
$${\mathbb H}_n=\,\big\{\,\Om\in \BC^{(n,n)}\,\big|\ \Om=\,^t\Om,\ \ \ \text{Im}\,\Om>0\,\big\}$$
be the Siegel upper half plane of degree $n$ and let
$$Sp(n,\BR)=\big\{ g\in \BR^{(2n,2n)}\ \big| \ ^t\!gJ_ng= J_n\ \big\}$$
be the symplectic group of degree $n$, where $F^{(k,l)}$ denotes
the set of all $k\times l$ matrices with entries in a commutative
ring $F$ for two positive integers $k$ and $l$, $^t\!M$ denotes
the transpose of a matrix $M,\ \text{Im}\,\Om$ denotes the
imaginary part of $\Om$ and
$$J_n=\begin{pmatrix} 0&I_n \\
                   -I_n&0 \\ \end{pmatrix}.$$
Here $I_n$ denotes the identity matrix of degree $n$.
We see that $Sp(n,\BR)$ acts on $\BH_n$ transitively by
\begin{equation*}g\cdot \Om=(A\Om+B)(C\Om+D)^{-1}, \end{equation*}
where $g=\begin{pmatrix} A&B\\ C&D\end{pmatrix}\in Sp(n,\BR)$ and
$\Om\in \BH_n.$

For two positive integers $n$ and $m$, we consider the Heisenberg
group
$$H_{\BR}^{(n,m)}=\{\,(\l,\mu;\k)\,|\ \l,\mu\in \BR^{(m,n)},\ \k\in \BR^{(m,m)},\ \
\k+\mu\,^t\l\ \text{symmetric}\ \}$$ endowed with the following
multiplication law
$$(\l,\mu;\k)\circ (\l',\mu';\k')=(\l+\l',\mu+\mu';\k+\k'+\l\,^t\mu'-
\mu\,^t\l').$$ We let
$$G^J=Sp(n,\BR)\ltimes H_{\BR}^{(n,m)}\quad \ ( \textrm{semi-direct product})$$
be the Jacobi group endowed with the following multiplication law
$$\Big(g,(\lambda,\mu;\kappa)\Big)\cdot\Big(g',(\lambda',\mu';\kappa')\Big) =\,
\Big(gg',(\widetilde{\lambda}+\lambda',\widetilde{\mu}+ \mu';
\kappa+\kappa'+\widetilde{\lambda}\,^t\!\mu'
-\widetilde{\mu}\,^t\!\lambda')\Big)$$ with $g,g'\in Sp(n,\BR),
(\lambda,\mu;\kappa),\,(\lambda',\mu';\kappa') \in
H_{\BR}^{(n,m)}$ and
$(\widetilde{\lambda},\widetilde{\mu})=(\lambda,\mu)g'$.
Then we have the {\it natural transitive action} of $G^J$ on the Siegel-Jacobi space $\BH_{n,m}:=\BH_n\times \BC^{(m,n)}$ defined by
\begin{equation*}
\Big(g,(\lambda,\mu;\kappa)\Big)\cdot (\Om,Z)=\Big((A\Om +B)(C\Om +D)^{-1},(Z+\lambda \,\Om+\mu)
(C\,\Om+D)^{-1}\Big),
\end{equation*}
where $g=\begin{pmatrix} A&B\\
C&D\end{pmatrix} \in Sp(n,\BR),\ (\lambda,\mu; \kappa)\in
H_{\BR}^{(n,m)}$ and $(\Om,Z)\in \BH_{n,m}.$ Thus $\BH_{n,m}$ is a homogeneous K{\"a}hler space which is not symmetric. In fact, $\BH_{n,m}$ is biholomorphic to
the homogeneous space $G^J/K^J$, where $K^J\cong U(n)\times S(m,\BR).$ Here $U(n)$ denotes the unitary group of degree $n$ and
$S(m,\BR)$ denote the abelian additive group consisting of all $m\times m$ symmetric real matrices.
We refer to \cite{BS, EZ, IOY},\,\cite{YJ6}-\cite{Zi} for more details on materials related to the Siegel-Jacobi space, e.g., Jacobi forms, invariant metrics,
invariant differential operators and Maass-Jacobi forms.

\vskip 0.2cm The Weil representation for a symplectic group was
first introduced by A. Weil in \cite{W} to reformulate Siegel's
analytic theory of quadratic forms (cf.\,\cite{Si}) in terms of
the group theoretical theory. It is well known that the Weil
representation plays a central role in the study of the
transformation behaviors of theta series.
In \cite{YJ17}, Yang constructed the Schr{\"o}dinger-Weil representation $\om_\CM$ of the Jacobi group $G^J$
associated with a positive definite symmetric real matrix $\CM$ of degree $n$ explicitly.

\vskip 0.21cm
This paper is organized as follows.
In Section 2, we review the Schr{\"o}dinger-Weil
representation $\omega_\CM$ of the Jacobi group $G^J$ associated
with a symmetric positive definite matrix $\CM$ and recall
the basic actions of $\omega_\CM$ on the representation space
$L^2\big(\BR^{(m,n)}\big)$ which were expressed explicitly in \cite{YJ17}. In Section 3, we define the theta sum $\Theta_f^{[\CM]}(\tau,\phi\,; \la,\mu,\kappa)$
of higher index and obtain some properties of the theta sum. The theta sum $\Theta_f^{[\CM]}(\tau,\phi\,; \la,\mu,\kappa)$ is a generalization of
the theta sum defined by J. Marklof \cite{Ma}.

\vskip 0.2cm \noindent {\bf Notations\,:} \ \ We denote by $\BZ,\,\,\BR$
and $\BC$ the ring of integers, the field of real numbers and the field of complex numbers
respectively. $\BC^{\times}$ denotes the multiplicative group of
nonzero complex numbers and $\BZ^{\times}$ denotes the set of all nonzero integers.
$T$ denotes the multiplicative group of
complex numbers of modulus one. The symbol ``:='' means that the
expression on the right is the definition of that on the left. For
two positive integers $k$ and $l$, $F^{(k,l)}$ denotes the set of
all $k\times l$ matrices with entries in a commutative ring $F$.
For a square matrix $A\in F^{(k,k)}$ of degree $k$, $\sigma(A)$
denotes the trace of $A$. For any $M\in F^{(k,l)},\ {}^t\!M$ denotes
the transpose of a matrix $M$. $I_n$ denotes the identity matrix of
degree $n$. We put $i=\sqrt{-1}.$ For a positive integer $m$ we denote by $S(m,F)$ the additive group consisting of all $m\times m$
symmetric matrices with coefficients in a commutative ring $F$.

\end{section}

\begin{section}{{\bf The Schr{\"o}dinger-Weil Representation}}
\setcounter{equation}{0}

In this section we review the Schr{\"o}dinger-Weil representation of the Jacobi group $G^J$ (cf. \cite{YJ17}, Section 3).\\
\indent
Throughout this section we assume that $\CM$ is a positive definite
symmetric real $m\times m$ matrix.
We let
\begin{equation*}
 L=\left\{\,(0,\mu;\kappa)\in H_{\mathbb
{R}}^{(n,m)}\,\Big| \, \mu\in \mathbb {R}^{(m,n)},\
\kappa=\,^t\!\kappa\in \mathbb {R}^{(m,m)}\ \right\}.
\end{equation*}
\noindent be a commutative normal subgroup of $H_{\mathbb{R}}^{(n,m)}$ and $\chi_\CM:L\lrt \BC^x$ be the unitary character of $L$
defined by
$$\chi_\CM ((0,\mu;\kappa)):=e^{\pi\,i\,\s (\CM \kappa)},\ \ \ (0,\mu;\kappa)\in L.$$
The representation ${\mathscr W}_\CM ={\rm Ind}_L^{ H_{\mathbb{R}}^{(n,m)}}\chi_\CM $ induced by $\chi_\CM$ from $L$ is realized on
the Hilbert space $H(\chi_\CM)\cong L^2\left(\mathbb{R}^{(m,n)}, d\xi\right)$.
${\mathscr W}_\CM$ is irreducible (cf.\,\cite{YJ1}, Theorem 3) and is called the Schr{\"o}dinger representation ${\mathscr W}_\CM$ of
the Heisenberg group $\hrnm$ with the central character $\chi_\CM$. We refer to \cite{YJ1, YJ2, YJ3, YJ4, YJ5, YB} for more details on representations of the Heisenberg
group $H_\BR^{(n,m)}$ and their related topics.
Then ${\mathscr W}_\CM$ is expressed explicitly as
\begin{equation}
\left[ {\mathscr W}_\CM (h_0)f\right](\lambda)=e^{\pi\,
i\,\sigma\{\CM(\kappa_0+\mu_0\,^t\!\lambda_0+
2\lambda\,^t\!\mu_0)\}}\,f(\lambda+\lambda_0),
\end{equation}
\noindent where $h_0=(\lambda_0,\mu_0;\kappa_0)\in H$ and
$\lambda\in\BR^{(m,n)}.$ See Formula (2.4) in \cite{YJ17} for more detail on ${\mathscr W}_\CM$. We
note that the symplectic group $Sp(n,\BR)$ acts on $\hrnm$ by
conjugation inside $G^J$. For a fixed element $g\in Sp(n,\BR)$,
the irreducible unitary representation ${\mathscr W}_\CM^g$ of
$\hrnm$ defined by
\begin{equation}
{\mathscr W}_\CM^g(h)={\mathscr W}_\CM(ghg^{-1}),\quad h\in\hrnm
\end{equation}
has the property that

\begin{equation*}
{\mathscr W}_\CM^g((0,0;\k))={\mathscr W}_\CM((0,0;\k))=e^{\pi
i\,\s(\CM \k)}\,\textrm{Id}_{H(\chi_\CM)},\quad \k\in
S(m,\BR).
\end{equation*}
Here $\textrm{Id}_{H(\chi_\CM)}$ denotes the identity operator on
the Hilbert space $H(\chi_\CM).$ According to Stone-von Neumann
theorem, there exists a unitary operator $R_\CM(g)$ on
$H(\chi_\CM)$  with $R_\CM (I_{2n})=\textrm{Id}_{H(\chi_\CM)}$ such that
\begin{equation}
R_\CM(g){\mathscr W}_\CM(h)={\mathscr
W}_\CM^g(h) R_\CM(g)\qquad {\rm for\ all}\ h\in\hrnm.
\end{equation}
We observe that
$R_\CM(g)$ is determined uniquely up to a scalar of modulus one.

\vskip 0.35cm
From now on, for brevity, we put $G=Sp(n,\BR).$ According to
Schur's lemma, we have a map $c_\CM:G\times G\lrt T$ satisfying
the relation
\begin{equation}
R_\CM(g_1g_2)=c_\CM(g_1,g_2)R_\CM(g_1)R_\CM(g_2)\quad \textrm{for
all }\ g_1,g_2\in G.
\end{equation}
We recall that $T$ denotes the multiplicative group of complex numbers of modulus one.
Therefore $R_\CM$ is a projective representation of $G$ on
$H(\chi_\CM)$ and $c_\CM$ defines the cocycle class in $H^2(G,T).$
The cocycle $c_\CM$ yields the central extension $G_\CM$ of $G$ by
$T$. The group $G_\CM$ is a set $G\times T$ equipped with the
following multiplication

\begin{equation}
(g_1,t_1)\cdot (g_2,t_2)=\big(g_1g_2,t_1t_2\,
c_\CM(g_1,g_2)^{-1}\,\big),\quad g_1,g_2\in G,\ t_1,t_2\in T.
\end{equation}
We see immediately that the map ${\widetilde R}_\CM:G_\CM\lrt
GL(H(\chi_\CM))$ defined by

\begin{equation}
{\widetilde R}_\CM(g,t)=t\,R_\CM(g) \quad \textrm{for all}\
(g,t)\in G_\CM
\end{equation}
is a {\it true} representation of $G_\CM.$ As in Section 1.7 in
\cite{LV}, we can define the map $s_\CM:G\lrt T$ satisfying the
relation
\begin{equation*}
c_\CM(g_1,g_2)^2=s_\CM(g_1)^{-1}s_\CM(g_2)^{-1}s_\CM(g_1g_2)\quad
\textrm{for all}\ g_1,g_2\in G.
\end{equation*}
Thus we see that
\begin{equation}
G_{2,\CM}=\left\{\, (g,t)\in G_\CM\,|\ t^2=s_\CM(g)^{-1}\,\right\}
\end{equation}

\noindent is the metaplectic group associated with $\CM$ that is a
two-fold covering group of $G$. The restriction $R_{2,\CM}$ of
${\widetilde R}_\CM$ to $G_{2,\CM}$ is the $\textsf{Weil representation}$ of
$G$ associated with $\CM$.
\par
If we identify $h=(\lambda,\mu;\kappa)\in \hrnm$ (resp. $g\in Sp(n, \BR)$) with
$(I_{2n},(\lambda,\mu;\kappa))\in G^J$ (resp. $(g,(0,0;0))\in G^J),$
every element $\tilde g$ of $G^J$ can be written as $\tilde g =hg$ with $h\in \hrnm$
and $g\in Sp(n, \BR)$. In fact,
\begin{equation*}
(g,(\la,\mu;\kappa))=(I_{2n},((\la,\mu)g^{-1};\kappa))\,(g,(0,0;0))=((\la,\mu)g^{-1};\kappa)\cdot g.
\end{equation*}
Therefore we define the {\it projective} representation $\pi_\CM$ of the Jacobi group
$G^J$ with cocycle $c_\CM (g_1,g_2)$ by
\begin{equation}
\pi_\CM(hg)={\mathscr W}_\CM(h)\,R_\CM(g),\quad h\in\hrnm,\ g\in G.
\end{equation}
\par
We let
$$G_{\CM}^J\!=G_{\CM}\ltimes \hrnm$$
be the semidirect product of $G_{\CM}$ and $\hrnm$ with the multiplication law
\begin{eqnarray*}
&&\big( (g_1,t_1),(\la_1,\mu_1;\kappa_1)\big)\cdot \big( (g_2,t_2),(\la_2,\mu_2\,;\kappa_2)\big)\\
&=&\big( (g_1,t_1)(g_2,t_2),(\tilde\la+\la_2,\tilde\mu+\mu_2\,;\kappa_1+\kappa_2+ \tilde\la\,^t\!\mu_2-\tilde\mu\,^t\!\la_2)\big),
\end{eqnarray*}
where $(g_1,t_1), (g_2,t_2)\in G_{\CM},\ (\la_1,\mu_1;\kappa_1), (\la_2,\mu_2\,;\kappa_2)\in\hrnm$ and
$(\tilde\la,\tilde\mu)=(\la,\mu)g_2.$
If we identify $h=(\lambda,\mu\,;\kappa)\in \hrnm$ (resp. $(g,t)\in G_{\CM})$ with
$((I_{2n},1),(\lambda,\mu\,;\kappa))\in G^J_{\CM}$ (resp. $((g,t),(0,0;0))\in G^J_{\CM}),$
we see easily that every element $\big( (g,t),(\la,\mu\,;\kappa)\big)$ of $G_{\CM}^J$ can be expressed as
$$\big( (g,t),(\la,\mu\,;\kappa)\big)=\big( (I_{2n},1),((\la,\mu)g^{-1};\kappa)\big) \big( (g,t),(0,0;0)\big)=((\la,\mu)g^{-1};\kappa)(g,t).  $$
Now we can define the {\it true} representation $\widetilde\om_\CM$ of $G_{\CM}^J$ by
\begin{equation}
\widetilde\omega_\CM(h\!\cdot\!(g,t))=t\,\pi_\CM(hg)=t\, {\mathscr
W}_\CM(h)\,R_\CM(g),\quad h\in\hrnm,\ (g,t)\in G_{\CM}.
\end{equation}

\vskip0.2cm
We recall that the following matrices
\begin{eqnarray*}
t(b)&=&\begin{pmatrix} I_n& b\\
                   0& I_n\end{pmatrix}\ \textrm{with any}\
                   b=\,{}^tb\in \BR^{(n,n)},\\
g(\alpha)&=&\begin{pmatrix} {}^t\alpha & 0\\
                   0& \alpha^{-1}  \end{pmatrix}\ \textrm{with
                   any}\ \alpha\in GL(n,\BR),\\
\s_n&=&\begin{pmatrix} 0& -I_n\\
                   I_n&\ 0\end{pmatrix}
\end{eqnarray*}
\noindent generate the symplectic group $G=Sp(n,\BR)$
(cf.\,\cite[p.\,326]{F},\,\cite[p.\,210]{Mum}). Therefore the
following elements $h_t(\lambda,\mu\,;\kappa),\
t(b\,;t),\,g(\alpha\,;t)$ and $\s_{n\,;t}$ of $G_\CM\ltimes \hrnm$
defined by
\begin{eqnarray*}
&& h_t(\la,\mu\,;\kappa)=\big( (I_{2n},t),(\la,\mu;\kappa)\big)\
\textrm{with}\ t\in T,\ \la,\mu\in
\BR^{(m,n)}\ \textrm{and}\ \kappa\in\BR^{(m,m)} ,\\
&&t(b\,;t)=\big( (t(b),t),(0,0;0) \big)\ \textrm{with any}\
                   b=\,{}^tb\in \BR^{(n,n)},\ t\in T,\\
&& g(\alpha\,;t)=\left(
\big(g(\alpha),t),(0,0;0)\right)\
\textrm{with any}\ \alpha\in GL(n,\BR)\ \textrm{and}\ t\in T,\\
 &&\s_{n\,;\,t}=\left( (\s_n,t),(0,0;0)\right)\
 \textrm{with}\ t\in T
\end{eqnarray*}
generate the group $G_\CM\ltimes\hrnm.$ We can show that the
representation ${\widetilde \om}_\CM$ is realized on the
representation $H(\chi_\CM)=L^2\big(\rmn\big)$ as follows: for
each $f\in L^2\big(\rmn\big)$ and $x\in \rmn,$ the actions of
${\widetilde \om}_\CM$ on the generators are given by

\begin{eqnarray}
\left[ {\widetilde \om}_\CM
\big(h_t(\lambda,\mu\,;\kappa)\big)f\right](x)&=&\,t\,e^{\pi
i\,\s\{\CM(\kappa+\mu\,{}^t\!\lambda+2\,x\,{}^t\mu)\}}\,f(x+\lambda),\\
\left[ {\widetilde \om}_\CM\big(t(b\,;t)\big)f\right](x)&=& t\,e^{\pi i\,\s(\CM\, x\,b\,{}^tx)}f(x),\\
\left[ {\widetilde \om}_\CM\big(g(\alpha\,;t)\big)f\right](x)&=& t\,| \det \alpha|^{\frac m2}\,f(x\,{}^t\alpha),
\end{eqnarray}
\begin{equation}
\left[ {\widetilde \om}_\CM\big(\s_{n\,;\,t}\big)f\right](x)=t\,
\big( \det \CM\big)^{\frac n2}\,\int_{\rmn}f(y)\,e^{-2\,\pi i\,\s(\CM\,y\,{}^tx)}\,dy.
\end{equation}
\par
Let
$$G_{2,\CM}^J\!=G_{2,\CM}\ltimes \hrnm$$
be the semidirect product of $G_{2,\CM}$ and $\hrnm$. Then $G_{2,\CM}^J$ is a subgroup of $G_\CM^J$ which is a two-fold covering group of the Jacobi group $G^J.$
The restriction $\om_\CM$ of $\widetilde \om_\CM$ to $G_{2,\CM}^J$ is called the $\textsf{Schr{\"o}dinger-Weil}$
$\textsf{representation}$ of $G^J$ associated with $\CM$.

\vskip 0.2cm\noindent
\begin{remark}
In the case $n=m=1,\ \omega_\CM$ is dealt in \cite{BS} and \cite{Ma}.
\end{remark}

\begin{remark}
The Schr{\"o}dinger-Weil representation is applied usefully to the theory of Maass-Jacobi forms \cite{P}.
\end{remark}

\end{section}

\vskip 1cm

\newcommand\wg{\widetilde g}
\newcommand\mfm{{\mathscr F}^{(\CM)} }
\newcommand\mfoz{{\mathscr F}^{(\CM)}_{\Om,Z} }
\newcommand\wgm{{\widetilde \gamma} }
\newcommand\Tm{\Theta^{(\CM)} }

\begin{section}{{\bf Theta Sums of Higher Index}}
\setcounter{equation}{0}

\vskip 0.35cm
Let $\mathcal M$ be a positive definite symmetric real matrix of degree $m$. We recall the Schr{\"o}dinger representation ${\mathscr W}_\CM$ of the Heisenberg
group $\hrnm$ associated with $\CM$ that is given by Formula (2.1) in Section 2. We note that for an element $(\la,\mu;\kappa)$ of $\hrnm$, we have the
decomposition
\begin{equation*}
(\la,\mu;\kappa)=(\la,0;0)\circ (0,\mu;0)\circ (0,0;\kappa\!-\!\la\, {}^t\!\mu).
\end{equation*}

We consider the embedding $\Phi_n :SL(2,\BR)\lrt Sp(n,\BR)$ defined by
\begin{equation}
\Phi_n \left( \begin{pmatrix} a & b\\ c & d \end{pmatrix}\right):=
\begin{pmatrix} aI_n & bI_n\\ cI_n & d I_n\end{pmatrix},\qquad \begin{pmatrix} a & b\\ c & d \end{pmatrix}\in SL(2,\BR).
\end{equation}

\vskip 0.35cm
For $x,y\in \BR^{(m,n)},$ we put
$$ (x,y)_\CM:=\s (\,{}^tx \CM y)\qquad \textrm{and}\qquad \| x\|_\CM :=\sqrt{(x,x)_\CM}.$$
According to Formulas (2.11)-(2.13), for any $ M=\begin{pmatrix} a & b\\ c & d \end{pmatrix}\in SL(2,\BR)\hookrightarrow Sp(n,\BR)$ and
$f\in L^2\left(\BR^{(m,n)}\right)$, we have the following explicit representation
\begin{equation}
[R_\CM (M)f](x)= \begin{cases} |a|^{\frac{mn}2} e^{ab \|x\|_\CM^2 \pi i} f(ax)  & \text{if $c=0$,}\\
(\det \CM)^{\frac n2}\, |c|^{-{\frac{mn}2}} \int_{\BR^{(m,n)}} e^{{\frac{\alpha(M,x,y,\CM)}c} \pi i} f(y) dy
& \text{if $c\neq 0$,}\end{cases}
\end{equation}
where
$$ \alpha(M,x,y,\CM)= a\,\|x\|_\CM^2 + d\, \|y\|_\CM^2 - 2 (x,y)_\CM.$$
Indeed, if $a=0$ and $c\neq 0$, using the decomposition
\begin{equation*}
M=\begin{pmatrix} 0 & -c^{-1}\\ c & d \end{pmatrix}=
\begin{pmatrix} 0 & -1\\ 1 & \ 0 \end{pmatrix} \begin{pmatrix} c & d\\ 0 & c^{-1} \end{pmatrix}
\end{equation*}
and
if $a\neq 0$ and $c\neq 0$, using the decomposition
\begin{equation*}
M=\begin{pmatrix} a & b\\ c & d \end{pmatrix}=\begin{pmatrix} a & c^{-1} \\ 0 & a^{-1} \end{pmatrix}
\begin{pmatrix} 0 & -1\\ 1 & \ 0 \end{pmatrix} \begin{pmatrix} ac & ad\\ 0 & (ac)^{-1} \end{pmatrix},
\end{equation*}
we obtain Formula (3.2).

\vskip 0.35cm
If
\begin{equation*}
M_1=\begin{pmatrix} a_1 & b_1\\ c_1 & d_1 \end{pmatrix},\quad M_2=\begin{pmatrix} a_2 & b_2\\ c_2 & d_2 \end{pmatrix}\quad
\textrm{and}\quad M_3=\begin{pmatrix} a_3 & b_3\\ c_3 & d_3 \end{pmatrix}\in SL(2,\BR)
\end{equation*}
with $M_3=M_1M_2$, the corresponding cocycle is given by
\begin{equation}
c_\CM (M_1,M_2)=e^{-i\, \pi\, mn\,\textrm{sign}(c_1c_2c_3)/4},
\end{equation}
where
\begin{equation*}
\textrm{sign}(x)= \begin{cases}  -1 \qquad &(x<0)\\
\ \ 0 \qquad &(x=0)\\
\ \ 1 \qquad &(x>0). \end{cases}
\end{equation*}
In the special case when
\begin{equation*}
M_1=\begin{pmatrix} \cos \phi_1 & -\sin \phi_1 \\ \sin \phi_1 & \ \ \cos \phi_1 \end{pmatrix}\quad \textrm{and}\quad
M_2=\begin{pmatrix} \cos \phi_2 & -\sin \phi_2 \\ \sin \phi_2 & \ \ \cos \phi_2 \end{pmatrix},
\end{equation*}
we find
\begin{equation*}
c_\CM (M_1,M_2)=e^{-i\, \pi\, mn\,(\sigma_{\phi_1}+\sigma_{\phi_2}-\sigma_{\phi_1+\phi_2})/4},
\end{equation*}
where
\begin{equation*}
\sigma_\phi= \begin{cases}  2\nu \qquad & \text{if $\phi=\nu\pi$}\\
2\nu+1 \qquad & \text{if $\nu\pi <\phi< (\nu+1)\pi.$} \end{cases}
\end{equation*}
It is well known that every $M\in SL(2,\BR)$ admits the unique Iwasawa decomposition
\begin{equation}
M=\begin{pmatrix} 1 & u\\ 0 & 1 \end{pmatrix} \begin{pmatrix} v^{1/2} & 0 \\ 0 & v^{-1/2} \end{pmatrix}
\begin{pmatrix} \cos \phi & -\sin \phi \\\sin \phi & \ \ \cos\phi \end{pmatrix},
\end{equation}
where $\tau=u+iv \in \BH_1$ and $\phi\in [0,2\pi).$ This parametrization $M=(\tau,\phi)$ in $SL(2,\BR)$ leads to the natural action of
$SL(2,\BR)$ on $\BH_1\times [0,2\pi)$ defined by
\begin{equation}
\begin{pmatrix} a & b\\ c & d \end{pmatrix}(\tau,\phi):=\left( \frac{a\tau+b}{c\tau+d},\, \phi + \textrm{arg} (c\tau+d)\ \textrm{mod}\ 2\pi \right).
\end{equation}

\begin{lemma}
For two elements $g_1$ and $g_2$ in $SL(2,\BR)$, we let
\begin{equation*}
g_1=\begin{pmatrix} 1 & u_1\\ 0 & 1 \end{pmatrix} \begin{pmatrix} v_1^{1/2} & 0 \\ 0 & v_1^{-1/2} \end{pmatrix}
\begin{pmatrix} \cos \phi_1 & -\sin \phi_1 \\\sin \phi_1 & \ \ \cos\phi_1 \end{pmatrix}
\end{equation*}
and
\begin{equation*}
g_2=\begin{pmatrix} 1 & u_2\\ 0 & 1 \end{pmatrix} \begin{pmatrix} v_2^{1/2} & 0 \\ 0 & v_2^{-1/2} \end{pmatrix}
\begin{pmatrix} \cos \phi_2 & -\sin \phi_2 \\\sin \phi_2 & \ \ \cos\phi_2 \end{pmatrix}
\end{equation*}
be the Iwasawa decompositions of $g_1$ and $g_2$ respectively, where $u_1,u_2\in\BR,\ v_1>0,
\,v_2>0$ and $0\leq \phi_1,\phi_2 < 2\pi.$ Let
\begin{equation*}
g_3=g_1g_2=\begin{pmatrix} 1 & u_3\\ 0 & 1 \end{pmatrix} \begin{pmatrix} v_3^{1/2} & 0 \\ 0 & v_3^{-1/2} \end{pmatrix}
\begin{pmatrix} \cos \phi_3 & -\sin \phi_3 \\\sin \phi_3 & \ \ \cos\phi_3 \end{pmatrix}
\end{equation*}
be the Iwasawa decomposition of $g_3=g_1g_2.$ Then we have
\begin{eqnarray*}
u_3&=& \frac{A}{(u_2\sin \phi_1+\cos \phi_1)^2+(v_2\sin\phi_1)^2},\\
v_3&=& \frac{v_1v_2}{(u_2\sin \phi_1+\cos \phi_1)^2+(v_2\sin\phi_1)^2}
\end{eqnarray*}
and
\begin{equation*}
\phi_3=tan^{-1} \left[
{ {(v_2\cos\phi_2+u_2 \sin\phi_2)\tan\phi_1 +\sin\phi_2}\over {(-v_2\sin\phi_2+u_2 \cos\phi_2)\tan\phi_1 +\cos\phi_2} }\right],
\end{equation*}
where
\begin{eqnarray*}
A&=&u_1 (u_2\sin \phi_1+\cos \phi_1)^2+ (u_1v_2-v_1u_2)\sin^2 \phi_1 \\
 && \ +\,v_1u_2\cos^2\phi_1+ v_1(u_2^2+v_2^2-1)\sin\phi_1\cos\phi_1.
\end{eqnarray*}
\end{lemma}

\vskip 0.26cm \noindent
{\it Proof.} If $g\in SL(2,\BR)$ has the unique Iwasawa decomposition (3.4), then we get the following
\begin{eqnarray*}
a&=& v^{1/2}\cos\phi +uv^{-1/2}\sin\phi,\\
b&=& -v^{1/2}\sin\phi +uv^{-1/2}\cos\phi,\\
c&=& v^{-1/2}\sin\phi, \quad   d=v^{-1/2}\cos\phi,\\
u&=&(ac+bd)\left(c^2+d^2\right)^{-1},\quad v=\left(c^2+d^2\right)^{-1},\quad \tan\phi={c\over d}\,    .
\end{eqnarray*}
We set
\begin{equation*}
g_3=g_1g_2=\begin{pmatrix} a_3 & b_3\\ c_3 & d_3 \end{pmatrix}.
\end{equation*}
Since
\begin{equation*}
u_3=(a_3c_3+b_3 d_3)\left(c_3^2+d_3^2\right)^{-1},\quad v=\left(c_3^2+d_3^2\right)^{-1},\quad \tan\phi_3={c_3\over d_3},
\end{equation*}
by an easy computation, we obtain the desired results.
\hfill $\square$

\vskip 0.53cm
Now we use the new coordinates $(\tau=u+iv,\phi)$ with $\tau\in\BH_1$ and $\phi\in [0,2\pi)$ in $SL(2,\BR).$ According to Formulas
(2.11)-(2.13), the projective representation $R_\CM$ of $SL(2,\BR)\hookrightarrow Sp(n,\BR)$ reads in these coordinates $(\tau=u+iv,\phi)$ as
follows:
\begin{equation}
\left[R_\CM (\tau,\phi)f\right](x)=v^{\frac{mn}4}\,e^{u \|x\|_\CM^2 \pi\, i} \left[R_\CM(i,\phi)f\right]\big(v^{1/2}x \big),
\end{equation}
where $f\in L^2\left( \BR^{(m,n)}\right),\ x\in \BR^{(m,n)}$ and
\begin{eqnarray}
 &\left[R_\CM(i,\phi)f\right](x) \hskip 9cm\nonumber\\
 =&\begin{cases}
 f(x) & \text{if $\phi\equiv 0$ mod $2\pi$,}\\
 f(-x) & \text{if $\phi\equiv \pi$ mod $2\pi$,}\\
 (\det\CM)^{\frac n2}\,|\sin\phi|^{-{{mn}\over 2}}\,\int_{\BR^{(m,n)}}e^{B(x,y,\phi,\CM)\pi i}\,f(y)dy
 & \text{if $\phi\not\equiv 0$ mod $\pi$}.
 \end{cases}
\end{eqnarray}
Here $$ B(x,y,\phi,\CM)= { {\left( \|x\|_\CM^2 +  \|y\|_\CM^2\right) \cos\phi - 2(x,y)_\CM} \over {\sin\phi}}. $$
Now we set
$$ S=\begin{pmatrix} 0 & -1\\ 1 & \ \ 0 \end{pmatrix}.$$
We note that
\begin{equation}
\left[ R_\CM \left( i, {\pi\over 2}\right)f\right](x)=\left[ R_\CM(S)f\right](x)=(\det\CM)^{\frac n2}\,\int_{\BR^{(m,n)}} f(y\,)\,e^{-2\, (x,\,y)_\CM\,\pi\,i}\,dy
\end{equation}
for $f\in L^2\left( \BR^{(m,n)}\right).$

\begin{remark}
For Schwartz functions $f\in \mathscr{S} \left(\BR^{(m,n)}\right),$ we have
\begin{equation*}
\lim_{\phi\lrt 0\pm} |\sin\phi|^{-{{mn}\over 2}}\, \int_{\BR^{(m,n)}}e^{B(x,y,\phi,\CM)\,\pi\, i}\,f(y)dy= e^{\pm i\,\pi\, mn/4}f(x)\neq f(x).
\end{equation*}
\noindent
Therefore the projective representation $R_\CM$ is not continuous at $\phi=\nu \pi\,(\nu\in\BZ)$ in general.
If we set
\begin{equation*}
\tilde{R}_\CM (\tau,\phi)= e^{-i\,\pi\, mn\sigma_\phi/4} R_\CM (\tau,\phi),
\end{equation*}
$\tilde{R}_\CM$ corresponds to a unitary representation of the double cover of $SL(2,\BR)$ (cf. Formula (2.6) and \cite{LV}).
This means in particular that
\begin{equation*}
\tilde{R}_\CM (i,\phi)\tilde{R}_\CM (i,\phi')=\tilde{R}_\CM (i,\phi+\phi'),
\end{equation*}
where $\phi\in [0,4\pi)$ parametrises the double cover of $SO(2)\subset SL(2,\BR).$
\end{remark}

\vskip 0.53cm
We observe that for any element $(g,(\la,\mu;\kappa))\in G^J$ with $g\in Sp(n,\BR)$ and $(\la,\mu;\kappa)\in \hrnm$, we have the following decomposition
\begin{equation*}
(g,(\la,\mu;\kappa))=(I_{2n},((\la,\mu)g^{-1};\kappa))\,(g,(0,0;0))=((\la,\mu)g^{-1};\kappa)\cdot g.
\end{equation*}
Thus $Sp(n,\BR)$ acts on $\hrnm$ naturally by
\begin{equation*}
g\cdot (\la,\mu;\kappa)=\left( (\la,\mu)g^{-1};\kappa\right),\qquad g\in Sp(n,\BR), \ (\la,\mu;\kappa)\in \hrnm.
\end{equation*}

\begin{definition}
For any Schwartz function $f\in \mathscr{S} \left(\BR^{(m,n)}\right),$ we define the function $\Theta_f^{[\CM]}$ on the Jacobi group
$SL(2,\BR)\ltimes \hrnm\hookrightarrow G^J$ by
\begin{equation}
\Theta_f^{[\CM]}(\tau,\phi\,;\la,\mu,\kappa):=\sum_{\om\in\BZ^{(m,n)}} \left[ \pi_\CM \left( (\la,\mu;\kappa)(\tau,\phi)\right)f\right] (\omega),
\end{equation}
where $(\tau,\phi)\in SL(2,\BR)$ and $(\la,\mu\,;\kappa)\in \hrnm$. The function $\Theta_f^{[\CM]}$ is called the theta sum of index $\CM$ associated
to a Schwartz function $f$. The projective representation $\pi_\CM$ of the Jacobi group $G^J$ was already defined by Formula (2.8).
More precisely, for $\tau=u+iv\in\BH_1$ and $(\la,\mu;\kappa)\in \hrnm$,
we have
\begin{eqnarray*}
&&\Theta_f^{[\CM]}(\tau,\phi\,;\la,\mu,\kappa)= v^{\frac{mn}4}\,\,e^{2\,\pi\,i\,\s(\CM(\kappa+\mu {}^t\la))}\\
&& \quad\times \sum_{\om\in\BZ^{(m,n)}}\, e^{\pi\,i\,\left\{ u \|\om+\la\|_\CM^2\,+\,2 (\om,\,\mu)_\CM \right\}}\,
\left[  R_\CM (i,\phi)f\right] \left( v^{1/2}(\omega+\l)\right).
\end{eqnarray*}
\end{definition}

\begin{lemma}
We set $f_\phi:=\tilde{R}_\CM (i,\phi)f$ for $f\in \mathscr{S} \left(\BR^{(m,n)}\right)$. Then for any $R>1$, there exists a constant $C_R$ such that
for all $x\in \rmn$ and $\phi\in\BR,$
$$ |f_\phi(x)| \leq C_R \,\left( 1+ \|x\|_\CM\right)^{-R}.$$
\end{lemma}

\vskip 0.25cm\noindent
{\it Proof.} Following the arguments in the proof of Lemma 4.3 in \cite{Ma}, pp.\,428-429, we get the desired result.
\hfill $\square$

\begin{theorem}[Jacobi 1]
Let $\CM$ be a positive definite symmetric integral matrix of degree $m$ such that $\CM \BZ^{(m,n)}=\BZ^{(m,n)}.$ Then
for any Schwartz function $f\in \mathscr{S} \left(\BR^{(m,n)}\right),$ we have
$$\Theta_f^{[\CM]}\left( -{1\over {\tau}}, \,\phi+\textrm{arg}\,\tau\,;-\mu,\la,\kappa \right)=\big(\det\CM\big)^{-{\frac n2}}\,  c_\CM(S,(\tau,\phi)) \,
\Theta_f^{[\CM]}(\tau,\phi\,;\la,\mu,\kappa),$$
where
$$c_\CM(S,(\tau,\phi)):=e^{i\,\pi mn \,\textrm{sign}(\sin\phi\,\sin (\phi+\arg \tau))}.$$
\end{theorem}

\vskip 0.251cm\noindent
{\it Proof.} First we recall that for any Schwartz function $\varphi\in \mathscr{S} \left(\BR^{(m,n)}\right),$ the Fourier transform $\mathscr F\varphi$ of $\varphi$
is given by
\begin{equation*}
\big( \mathscr F \varphi\big)(x)=\int_{\rmn} \varphi(y)\, e^{-2\pi i\,\s ( y\,^t\!x)} dy.
\end{equation*}
Now we put
\begin{equation*}
 S=\begin{pmatrix} 0 & -1 \\ 1 & \ \ 0 \end{pmatrix}\in SL(2,\BZ) \hookrightarrow Sp(n,\BR)
\end{equation*}
and for any $F\in \mathscr{S} \left(\BR^{(m,n)}\right),$ we put
\begin{equation*}
F_\CM(x):=F(\CM^{-1}x), \quad x\in \rmn.
\end{equation*}
According to Formula (2.13), for any $F\in \mathscr{S} \left(\BR^{(m,n)}\right),$
\begin{eqnarray*}
\left[ R_\CM (S)F\right](x)&=& \big(\det\CM \big)^{\frac n2} \,\int_{\rmn} F(y)\,e^{-2\pi i\,\s (\CM y\,^t\!x)} dy \\
&=& \big(\det\CM \big)^{-\frac n2}\,\int_{\rmn} F(\CM^{-1}y)\,e^{-2\pi i\,\s ( y\,^t\!x)} dy \\
&=& \big(\det\CM \big)^{-\frac n2}\, \int_{\rmn} F_\CM (y)\,e^{-2\pi i\,\s ( y\,^t\!x)} dy\\
&=& \big(\det\CM \big)^{-\frac n2}\, \left[\mathscr F F_\CM\right](x).
\end{eqnarray*}
Thus we have
\begin{equation}
\mathscr F F_\CM=\big(\det\CM \big)^{\frac n2}\,R_\CM (S)F  \qquad \textrm{for}\  F\in \mathscr{S} \left(\BR^{(m,n)}\right).
\end{equation}
By Lemma 3.1, we get easily
\begin{equation}
S\cdot (\tau,\phi)=\left( -{1\over\tau}, \phi+ \arg \tau \right).
\end{equation}
If we take $F=\pi_\CM ((\la,\mu\,;\kappa)(\tau,\phi))f$ for $f\in \mathscr{S} \left(\BR^{(m,n)}\right)$, a fixed element
$(\la,\mu\,;\kappa)\in\hrnm$ and an fixed element $(\tau,\phi)\in SL(2,\BR),$ then it is easily seen that $F\in \mathscr{S} \left(\BR^{(m,n)}\right)$.
\par
According to Formulas (3.11), if  we take $F=\pi_\CM ((\la,\mu\,;\kappa)(\tau,\phi))f$ for $f\in \mathscr{S} \left(\BR^{(m,n)}\right)$,
\begin{eqnarray*}
\big[ R_\CM (S)F\big](x)
&=& \left[ R_\CM (S) \pi_\CM \big( (\la,\mu\,;\kappa)(\tau,\phi)\big)f\right](x),\quad x\in \rmn\\
&=& \left[ R_\CM (S) \mathscr W_\CM (\la,\mu\,;\kappa) R_\CM(\tau,\phi)f\right](x)\\
&=& \left[ \mathscr W_\CM \big( (\la,\mu)S^{-1};\kappa\big) R_\CM(S) R_\CM(\tau,\phi)f\right](x)\\
&=& c_\CM(S,(\tau,\phi))^{-1}\,\left[ \mathscr W_\CM  (-\mu,\la\,;\kappa) R_\CM\big( S\cdot (\tau,\phi)\big)f\right](x)\\
&=& c_\CM(S,(\tau,\phi))^{-1}\,\left[ \mathscr W_\CM  (-\mu,\la\,;\kappa) R_\CM \left( -{1\over\tau}, \phi+ \arg \tau \right)f\right](x)\\
&=& c_\CM(S,(\tau,\phi))^{-1}\,\left[ \pi_\CM \left( (-\mu,\la\,;\kappa)  \left( -{1\over\tau}, \phi+ \arg \tau \right)\right)f\right](x).
\end{eqnarray*}
Thus we obtain
\begin{equation}
\big[ R_\CM (S)F\big](x)= c_\CM(S,(\tau,\phi))^{-1}\,\left[ \pi_\CM \left( (-\mu,\la\,;\kappa)  \left( -{1\over\tau}, \phi+ \arg \tau \right)\right)f\right](x) .
\end{equation}
According to Poisson summation formula, we have
\begin{equation}
\sum_{\om\in\BZ^{(m,n)}}\left[ \mathscr F F_\CM \right](\om)=\sum_{\om\in\BZ^{(m,n)}} F_\CM (\om).
\end{equation}
It follows from (3.10) and (3.12) that
\begin{eqnarray*}
\sum_{\om\in\BZ^{(m,n)}}\left[ \mathscr F F_\CM \right](\om) &=& \big(\det\CM \big)^{\frac n2} \,\sum_{\om\in\BZ^{(m,n)}}\big[ R_\CM (S)F\big](\om)\\
&=& \big(\det\CM \big)^{\frac n2} \,c_\CM(S,(\tau,\phi))^{-1}\,\\
& & \times \sum_{\om\in\BZ^{(m,n)}} \,\left[ \pi_\CM \left( (-\mu,\la\,;\kappa)  \left( -{1\over\tau}, \phi+ \arg \tau \right)\right)f\right](x)\\
&=& \big(\det\CM \big)^{\frac n2} \,  c_\CM(S,(\tau,\phi))^{-1}\, \Theta_f^{[\CM]}\left( -{1\over {\tau}}, \,\phi+\textrm{arg}\,\tau\,;-\mu,\la,\kappa \right).
\end{eqnarray*}
On the other hand,
\begin{eqnarray*}
\sum_{\om\in\BZ^{(m,n)}} F_\CM (\om) &=& \sum_{\om\in\BZ^{(m,n)}}  F\big(\CM^{-1}\om \big) \\
&=& \sum_{\om\in\BZ^{(m,n)}} \left[ \pi_\CM ((\la,\mu\,;\kappa)(\tau,\phi))f \right] \big(\CM^{-1}\om \big)\\
&=& \sum_{\om\in\BZ^{(m,n)}} \left[ \pi_\CM ((\la,\mu\,;\kappa)(\tau,\phi))f \right] (\om)\quad \left( \because \ \CM^{-1}\BZ^{(m,n)}=\BZ^{(m,n)} \right)\\
&=& \Theta_f^{[\CM]}(\tau,\phi\,;\la,\mu,\kappa).
\end{eqnarray*}
Hence from (3.13) we obtain the desired formula
\begin{equation*}
\Theta_f^{[\CM]}\left( -{1\over {\tau}}, \,\phi+\textrm{arg}\,\tau\,;-\mu,\la,\kappa \right)= \big(\det\CM \big)^{-\frac n2} \, c_\CM(S,(\tau,\phi)) \,
\Theta_f^{[\CM]}(\tau,\phi\,;\la,\mu,\kappa).
\end{equation*}
If
$$ S=\begin{pmatrix} a_1 & b_1 \\ c_1 & d_1 \end{pmatrix},\quad
(\tau,\phi)=\begin{pmatrix} a_2 & b_2 \\ c_2 & d_2 \end{pmatrix}\quad \textrm{and}\quad
S\cdot (\tau,\phi)=\begin{pmatrix} a_3 & b_3 \\ c_3 & d_3 \end{pmatrix}\in SL(2,\BR),$$
according to Lemma 3.1, we get easily
$$c_1c_2c_3 = \big( u^2+v^2\big)^{1/2} \sin\phi\,\sin (\phi+\arg \tau),$$
where
$$ (\tau,\phi)=\begin{pmatrix} 1 & u\\ 0 & 1 \end{pmatrix} \begin{pmatrix} v^{1/2} & 0 \\ 0 & v^{-1/2} \end{pmatrix}
\begin{pmatrix} \cos \phi & -\sin \phi \\\sin \phi & \ \ \cos\phi \end{pmatrix}$$
is the Iwasawa decomposition of $(\tau,\phi)\in SL(2,\BR).$ Thus we obtain
$$c_\CM(S,(\tau,\phi))=e^{i\,\pi mn \,\textrm{sign}(c_1c_2c_3)}=     e^{i\,\pi mn \,\textrm{sign}(\sin\phi\,\sin (\phi+\arg \tau))}.$$
This completes the proof.
\hfill $\square$

\vskip 0.53cm
\begin{theorem}[Jacobi 2]
Let $\CM=(\CM_{kl})$ be a positive definite symmetric integral $m\times m$ matrix and
let $s=(s_{kj})\in \BZ^{(m,n)}$ be integral. Then we have
\begin{equation*}
\Theta_f^{[\CM]}(\tau+2,\phi\,;\la,s-2\,\la+\mu,\kappa-s\,^t\la)=\Theta_f^{[\CM]}(\tau,\phi\,;\la,\mu,\kappa)
\end{equation*}
for all $(\tau,\phi)\in SL(2,\BR)$ and $(\la,\mu;\kappa)\in \hrnm$.
\end{theorem}
\vskip 0.251cm\noindent
{\it Proof.} For brevity, we put $T_*=\begin{pmatrix} 1 & 2\\ 0 & 1 \end{pmatrix}$. According to Lemma 3.1, for any $(\tau,\phi)\in SL(2,\BR),$
the multiplication of $T_*$ and $(\tau,\phi)$ is given by
\begin{equation}
T_* (\tau,\phi)=(\tau+2,\phi).
\end{equation}
For $s\in \rmn,\ (\la,\mu\,;\kappa)\in\hrnm$ and $(\tau,\phi)\in SL(2,\BR),$ according to (3.14),
\begin{eqnarray*}
&& \pi_\CM ((0,s;0)T_*)\, \pi_\CM ((\la,\mu\,;\kappa)(\tau,\phi))\\
&=& \mathscr W_\CM (0,s;0)R_\CM (T_*)  \mathscr W_\CM (\la,\mu\,;\kappa) R_\CM (\tau,\phi)\\
&=& \mathscr W_\CM (0,s;0) \mathscr W_\CM \big((\la,\mu)T_*^{-1}\,;\kappa\big) R_\CM (T_*)R_\CM (\tau,\phi)\\
&=& c_\CM (T_*, (\tau,\phi))^{-1}  \mathscr W_\CM (\la, s-2\,\la+\mu\,;\kappa-s\,{}^t\!\la) R_\CM \big(T_* (\tau,\phi)\big)\\
&=&  \mathscr W_\CM (\la,s-2\,\la+\mu\,;\kappa-s\,{}^t\!\la) R_\CM (\tau+2,\phi)\\
&=& \pi_\CM \big( (\la,s-2\,\la+\mu\,;\kappa-s\,{}^t\!\la) (\tau+2,\phi) \big).
\end{eqnarray*}
Here we used the fact that $c_\CM (T_*, (\tau,\phi))=1$ because $T_*$ is upper triangular.
\par
On the other hand,
according to the assumptions on $\CM$ and $s$, for $f\in \mathscr{S} \left(\BR^{(m,n)}\right)$ and $\om\in \BZ^{(m,n)},$
using Formulas (2.1), (2.11) or (3.6), we have
\begin{eqnarray*}
&&\left[ \pi_\CM \big( (0,s;0)T_*\big) \,\pi_\CM \big( (\la,\mu;\kappa)(\tau,\phi)\big)f\right](\om)\\
&=& \left[ \mathscr W_\CM (0,s;0) R_\CM (T_*) \,\pi_\CM \big( (\la,\mu;\kappa)(\tau,\phi)\big)f\right](\om)\\
&=& e^{2\pi i \,\s(\CM \om\,{}^t\!s)}\cdot e^{2\,\|\om\|_\CM^2 \pi\, i}\,\left[ R_\CM(i,0)\,\pi_\CM \big( (\la,\mu;\kappa)(\tau,\phi)\big)f\right](\om)\\
&=& \left[ \pi_\CM \big( (\la,\mu;\kappa)(\tau,\phi)\big)f\right](\om).
\end{eqnarray*}
Here we used the facts that
$$ e^{2\pi i\, \s(\CM \om\,{}^t\!s)}=1,\quad   e^{2\,\|\om\|_\CM^2 \pi\, i}=1 \quad \textrm{and}\quad R_\CM(i,0)f=f\ (\textrm{cf}.\ (3.7)).$$
Therefore for $f\in \mathscr{S} \left(\BR^{(m,n)}\right)$,
\begin{eqnarray*}
&& \Theta_f^{[\CM]}(\tau+2,\phi\,;\la,s-2\,\la+\mu,\kappa-s\,^t\la)\\
&=&\sum_{\om\in\BZ^{(m,n)}} \left[ \pi_\CM \big( (\la,s-2\,\la+\mu,\kappa-s\,^t\la)(\tau+2,\phi)\big) f \right](\om)\\
&=& \sum_{\om\in\BZ^{(m,n)}} \left[ \pi_\CM \big( (0,s;0)T_*\big) \,\pi_\CM \big( (\la,\mu;\kappa)(\tau,\phi)\big)f\right](\om)\\
&=&\sum_{\om\in\BZ^{(m,n)}} \left[ \pi_\CM \big( (\la,\mu;\kappa)(\tau,\phi)\big)f\right](\om)\\
&=&\Theta_f^{[\CM]}(\tau,\phi\,;\la,\mu,\kappa).
\end{eqnarray*}
This completes the proof.
\hfill $\square$

\vskip 0.53cm
\begin{theorem}[Jacobi 3]
Let $\CM=(\CM_{kl})$ be a positive definite symmetric integral $m\times m$ matrix and let $(\la_0,\mu_0;\kappa_0)\in H_\BZ^{(m,n)}$
be an integral element of $\hrnm.$ Then we have
\begin{eqnarray*}
&&\Theta_f^{[\CM]}(\tau,\phi\,;\la+\la_0,\mu+\mu_0,\kappa+\kappa_0+\la_0\,{}^t\mu-\mu_0\,^t\la)\\
&=& e^{\pi\,i\,\s(\CM (\kappa_0+\mu_0\,{}^t\la_0))}
\Theta_f^{[\CM]}(\tau,\phi\,;\la,\mu,\kappa)
\end{eqnarray*}
for all $(\tau,\phi)\in SL(2,\BR)$ and $(\la,\mu;\kappa)\in \hrnm$.
\end{theorem}
\vskip 0.251cm\noindent
{\it Proof.} For any $f\in \mathscr{S} \left(\BR^{(m,n)}\right)$, we have
\begin{eqnarray*}
&& \sum_{\om\in\BZ^{(m,n)}} \left[ \mathscr W_\CM (\la_0,\mu_0;\kappa_0)\pi_\CM \big( (\la,\mu\,;\kappa)(\tau,\phi)\big)f\right](\om)\\
&=& \sum_{\om\in\BZ^{(m,n)}} \left[ \mathscr W_\CM (\la_0,\mu_0;\kappa_0)\mathscr W_\CM  (\la,\mu\,;\kappa) R_\CM(\tau,\phi) f\right](\om)\\
&=& \sum_{\om\in\BZ^{(m,n)}} \left[ \mathscr W_\CM (\la_0+\la,\mu_0+\mu;\kappa_0+\kappa+ \la_0\,{}^t\!\mu-\mu_0\,{}^t\!\la)) R_\CM(\tau,\phi) f\right](\om)\\
&=& \sum_{\om\in\BZ^{(m,n)}} \left[ \pi_\CM\big( (\la_0+\la,\mu_0+\mu;\kappa_0+\kappa+ \la_0\,{}^t\!\mu-\mu_0\,{}^t\!\la)(\tau,\phi)\big) f\right](\om)\\
&=& \Theta_f^{[\CM]}(\tau,\phi\,;\la+\la_0,\mu+\mu_0,\kappa+\kappa_0+\la_0\,{}^t\mu-\mu_0\,^t\la).
\end{eqnarray*}
On the other hand, for any $f\in \mathscr{S} \left(\BR^{(m,n)}\right)$, we have
\begin{eqnarray*}
&& \sum_{\om\in\BZ^{(m,n)}} \left[ \mathscr W_\CM (\la_0,\mu_0;\kappa_0)\pi_\CM \big( (\la,\mu\,;\kappa)(\tau,\phi)\big)f\right](\om)\\
&=& \sum_{\om\in\BZ^{(m,n)}} e^{\pi i\sigma\{\CM (\kappa_0+\mu_0\,^t\!\lambda_0+ 2\,\om \,^t\!\mu_0)\}} \left[ \pi_\CM (\tau,\phi\,; \la,\mu,\kappa)f\right](\om+\la_0)\\
&=& e^{\pi i\sigma\{\CM (\kappa_0+\mu_0\,^t\!\lambda_0\}}  \sum_{\om\in\BZ^{(m,n)}}  \left[ \pi_\CM (\tau,\phi\,; \la,\mu,\kappa)f\right](\om+\la_0)\quad
(\because \ \mu_0\ \textrm{is\ integral})\\
&=& e^{\pi i\sigma\{\CM (\kappa_0+\mu_0\,^t\!\lambda_0\}}  \sum_{\om\in\BZ^{(m,n)}}  \left[ \pi_\CM (\tau,\phi\,; \la,\mu,\kappa)f\right](\om)\quad
(\because \ \la_0\ \textrm{is\ integral})\\
&=&  e^{\pi i\sigma\{\CM (\kappa_0+\mu_0\,^t\!\lambda_0\}}\,\Theta_f^{[\CM]}(\tau,\phi\,;\la,\mu,\kappa).
\end{eqnarray*}
Finally we obtain the desired result.
\hfill $\square$

\vskip 0.35cm
We put $V(m,n)=\BR^{(m,n)}\times \BR^{(m,n)}$. Let
\begin{equation*}
G^{(m,n)}:=SL(2,\BR) \ltimes V(m,n)
\end{equation*}
be the group with the following multiplication law
\begin{equation}
(g_1,(\la_1,\mu_1))\cdot (g_2,(\la_2,\mu_2))=(g_1g_2,(\la_1,\mu_1)g_2+ (\la_2,\mu_2)),
\end{equation}
where $g_1,g_2\in SL(2,\BR)$ and $\la_1,\la_2,\mu_1,\mu_2\in \BR^{(m,n)}$.

\vskip 0.25cm\noindent
We define
$$  \G^{(m,n)}:= SL(2,\BZ)\ltimes H_\BZ^{(n,m)}.$$
Then $\G^{(m,n)}$ acts on $G^{(m,n)}$ naturally through the multiplication law (3.15).
\begin{lemma}
$\G^{(m,n)}$ is generated by the elements
$$ ( S, (0,0)),\quad (T_\flat,(0,s)) \quad \textrm{and}\quad (I_2,(\la_0,\mu_0)),$$
where
$$  S=\begin{pmatrix} 0 & -1 \\ 1 & \ \ 0 \end{pmatrix},\quad T_\flat=\begin{pmatrix} 1 & 1 \\ 0 & 1 \end{pmatrix}\quad \textrm{and}\  s,\la_0,\mu_0\in \BZ^{(m,n)}.$$
\end{lemma}
\vskip 0.251cm\noindent
{\it Proof.}
Since $SL(2,\BZ)$ is generated by $S$ and $T_\flat$, we get the desired result.
\hfill $\square$

\vskip 0.53cm
\noindent We define
\begin{eqnarray*}
&&\Theta_f^{[\CM]}(\tau,\phi;\la,\mu)\\
&=&  v^{\frac{mn}4}\,\sum_{\om\in\BZ^{(m,n)}}\, e^{\pi\,i\,\left\{ u \|\om+\la\|_\CM^2\,+\,2 (\om,\,\mu)_\CM \right\}}\,
\left[ R_\CM (i,\phi)f\right] \left( v^{1/2}(\omega+\l)\right).\nonumber
\end{eqnarray*}

\begin{theorem}
Let $\G^{(m,n)}_{[2]}$ be the subgroup of $\G^{(m,n)}$ generated by the elements
$$ ( S, (0,0)),\quad (T_*,(0,s)) \quad \textrm{and}\quad (I_2,(\la_0,\mu_0)),$$
where
$$  T_*=\begin{pmatrix} 1 & 2 \\ 0 & 1 \end{pmatrix}\quad \textrm{and}\  s,\la_0,\mu_0\in \BZ^{(m,n)}.$$
Let $\CM=(\CM_{kl})$ be a positive definite symmetric unimodular integral $m\times m$ matrix such that $\CM \BZ^{(m,n)}=\BZ^{(m,n)}.$
Then for $f,g\in \mathscr{S} \left( \BR^{(m,n)}\right),$ the function
\begin{equation*}
\Theta_f^{[\CM]}(\tau,\phi;\la,\mu)\,\overline{\Theta_g^{[\CM]}(\tau,\phi;\la,\mu)}
\end{equation*}
is invariant under the action of $\G^{(m,n)}_{[2]}$ on $G^{(m,n)}$.
\end{theorem}
\vskip 0.251cm\noindent
{\it Proof.} The proof follows directly from Theorem 3.1 (Jacobi 1), Theorem 3.2 (Jacobi 2) and Theorem 3.3 (Jacobi 3) because the left actions of the generators of
$\G^{(m,n)}_{[2]}$ are given by
\begin{eqnarray*}
&&((\tau,\phi),(\la,\mu))\longmapsto \left( \left( -{1\over \tau},\phi+\arg\tau\right),(-\mu,\la)\right),\\
&&((\tau,\phi),(\la,\mu)) \longmapsto ( (\tau+2,\phi),(\la,s-2\,\la+\mu))
\end{eqnarray*}
and
$$ ((\tau,\phi),(\la,\mu)) \longmapsto  ((\tau,\phi),(\la+\la_0,\mu+\mu_0)). $$
\hfill $\square$

\end{section}

\vskip 1cm
\bibliography{central}

\end{document}